\documentclass[a4paper,10pt]{article}
\usepackage[utf8x]{inputenc}
\usepackage{amssymb}

\title{Special Bohr - Sommerfeld geometry: variations}
\author{ Nikolay Tyurin \\
BLTPh JINR (Dubna) and MI RAS (Moscow)}

\date{}

\begin{document}

\maketitle

\begin{abstract} In the paper we continue to study Special Bohr - Sommerfeld geometry of compact symplectic manifolds. Using natural deformation
parameters we avoid the difficulties appeared in the definition of the moduli space of Special Bohr - Sommerfeld cycles for  compact simply connected algebraic varieties.
As a byproduct we present certain remarks on the Weinstein structures and Eliashberg conjectures.

\end{abstract}

\section*{Introduction}

Recall the basic constructions of Special Bohr - Sommerfeld geometry of compact symplectic manifolds (the details can be found in [1], [2], [3], [4]).

Let $(M, \omega)$ be a compact simply connected symplectic manifold of dimension $2n$ such that the symplectic form $\omega$ is of integer type so its cohomology class
is integer. Consider the corresponding complex line bundle $L \to M$ such that $c_1(L) = [\omega]$, equipped with a hermitian structure
$h$, so the corresponding space of hermitian connections ${\cal A}_h(L)$ contains a subset ${\cal O}(\omega)$ formed by the solutions
of the equation $F_a = 2 \pi i \omega$ (in the simply connected case the subset ${\cal O}(\omega)$ is an orbit of the gauge group action).
Choosing an element $a \in {\cal O}(\omega)$ one gets the corresponding prequantization pair $(L, a)$ which plays the key role in
Geometric Quantization procedure; from the GQ story we take the Hilbert space $\Gamma (M, L)$ consists of smooth sections of $L$
where the scalar product $<s_1, s_2> = \int_M (s_1, s_2)_h d \mu_L$ is generated by $h$ and the Lioville volume form $d \mu_L$.

Fixing a topological type $\rm{top} S$ of a smooth orientable  $n$ - dimensional manifold and a middle homology class $[S] \in H_n(M, \mathbb{Z})$
one gets the moduli space ${\cal B}_S$ of Bohr - Sommerfeld lagrangian cycles of fixed type, see [5], which is an infinite dimensional
Frechet smooth manifold, whose "points" can be understood as lagrangian submanifolds $S \subset M$ of fixed topological type which
satisfy the Bohr - Sommerfeld condition: for each $a \in {\cal O}(\omega)$ the pair $(L,a)$ admits covariantly constant sections being
restricted to $S$. The moduli space ${\cal B}_S$ subjects to another approach to Geometric Quantization called Lagrangian approach to GQ, see [7].
Here and below we are interested in smooth lagrangian submanifolds only.

Leaving aside GQ motivation and utilities, define a universal object ${\cal U}_{SBS}(a)$ in the direct product
$\mathbb{P} \Gamma (M, L) \times {\cal B}_S$ by the following rules. Pair $(p = [\alpha], S)$, where $\alpha \in \Gamma(M, L)$ is a section, representing
 the equivalence class $p$,  belongs to ${\cal U}_{SBS}(a)$ if the proportionality coefficient
$\frac{\alpha|_S}{\sigma_S}$ has the form $e^{i c} f$ where $\sigma_S$ is a covariantly
constant section of $(L,a)|_S$, $c$ is a real constant and $f \in C^{\infty}(S, \mathbb{R}_+)$ is a strictly positive real function on $S$.
Note that the changing of $a \in {\cal O}(\omega)$ is reflected in the definition since the covariantly constant sections $\sigma_S$ depend
on the choice of $a$ for any particular Bohr - Sommerfeld submanifold $S$.

By the very definition ${\cal U}_{SBS}(a)$ carries two natural projections:
$$
\mathbb{P} \Gamma (M, L) \leftarrow {\cal U}_{SBS}(a) \rightarrow {\cal B}_S,
$$
with the following properties. The first projection $p_1: {\cal U}_{SBS}(a) \to \mathbb{P} \Gamma (M, L)$ has discrete fibers;
its image is an open subset in $\mathbb{P}\Gamma (M, L)$; its differential never vanishes (here the smoothness of $S$ is crucial);
no ramification takes place in the picture (the detailed proofs in [jT1]). Since the projective space $\mathbb{P} \Gamma (M, L)$
carries the standard Kahler from $\Omega_{FS}$ of the Fubini - Study metric, it follows that in absolutely general situation
${\cal U}_{SBS}(a)$ is weakly Kahler manifold.

The second projection $p_2: {\cal U}_{SBS}(a) \to {\cal B}_S$ splits as $\pi \circ \tau$ where $\pi: T {\cal B}_S \to {\cal B}_S$ is the canonical projection
of the tangent bundle,  and
$\tau: {\cal U}_{SBS}(a) \to T {\cal B}_S$ is a map with Kahler fibers (the detailed proofs in [4]).

Now if we turn to a natural problem appeared many times both in symplectic geometry and mathematical physics which looks for
finite dimensional moduli spaces of lagrangian submanifolds satisfying certain additional conditions then one can see
that SBS geometry can be exploited in this way. Namely suppose one fixes a natural finite dimensional projective subspace $\mathbb{P}^N$
in $\mathbb{P} \Gamma (M, L)$ and then take its preimage under the first projection $p_1$ then such $p^{1}(\mathbb{P}^N)$
must be finite dimensional Kahler manifold. Moreover, the choice of an appropriate $\mathbb{P}^N \subset \mathbb{P} \Gamma (M, L)$
can be made almost automatically in the situation when $(M, \omega)$ is equipped with a compatible integrable complex structure $I$.
Indeed, since $M$ is compact then the holomorphic section space $H^0(M_I, L)$ of the prequantization bundle is finite dimensional,
therefore we can take $\mathbb{P}^N = \mathbb{P} H^0(M_I, L) \subset \mathbb{P} \Gamma (M, L)$ and then it is natural to define
 the preimage $p_1^{-1}(\mathbb{P} H^0(M_I, L))$ as a distinguished geometrical object which characterizes the Kahler nature
 of $(M, \omega, I)$.

 The construction is rather concise in the case of algebraic varieties. Indeed, a compact algebraic variety $X$
 with a very ample line bundle  $L$ (note that such  $L$ must exist by the very definition)  from the
 real geometry viewpoint is presented by $(M, \omega, I)$ where $c_1(L) = [\omega]$; this correspondence can be realized via the choice
 of an appropriate hermitian structure $h$ on $L$ such that $\omega = - \rm{d} I \rm{d} (\rm{ln} \vert \alpha \vert_h)$
 on the complement to the zeroset $D_{\alpha} = \{ \alpha = 0 \} \subset X$ for a holomorphic section $\alpha \in
 H^0(X, L)$. Clearly such $\omega$ is not unique, but the lagrangian geometries for different forms of this type are equivalent
 (at the same time for different principal polarizations $L_1$ and $L_2$  the corresponding pictures  can be different).
 In the presence of complex structure $I$ one has a distinguished hermitian connection $a_I \in {\cal O}(\omega)$,
 compatible with the holomorphic structure on $L$, and SBS - construction leads to the coarse definition:
 the moduli space ${\cal M}_{SBS}(c_1(L))$ of SBS lagrangian cycles is given by $p^{-1}(\mathbb{P} H^0(M_I, L)) \subset {\cal U}_{SBS}(a_I)$.

 The theory says  that this moduli space should be finite dimensional smooth Kahler variety, but this expectation is failed due to the following facts.
 In [2] one established that Bohr - Sommerfeld lagrangian submanifold $S$ is special with respect to a holomorphic section $\alpha \in H^0(M_I, L)$
 when one fixes the connection $a_I$ for the definition of ${\cal U}_{SBS}$
 if and only if it is contained by the Weinstein skeleton of the complement $M_I \backslash D_{\alpha}$
 where the skeleton is formed by finite trajectories of the gradient flow of function $- \rm{ln} \vert \alpha \vert_h$, while the zero divisor $D_{\alpha}$
 of section $\alpha$ attracts all infinite trajectories, see [6]. But as one claims no smooth lagrangian $S$ exists even in the simplest cases
 (in [3] we present the example when $M_I = \mathbb{C} \mathbb{P}^1$ and $L = {\cal O}(3)$, --- even there no smooth loops in the Weinstein skeleton
 for generic holomorphic section!). Therefore the coarse definition leads to the trivial answer.

 At the same time a parallel result was presented in [3]: for the same situation one constructs the moduli space of D - exact lagrangian submanifolds
 which is a Kahler manifold fibered over an open part of $\mathbb{P} H^0(M_I, L)$ with discrete fibers. This moduli space was denoted
 as $\tilde {\cal M}_{SBS} (c_1(L))$ since in [3] one claims that these moduli spaces  must be somehow related to Special Bohr - Sommerfeld geometry.

 The main aim of the present paper is to correct the coarse definition for the moduli space ${\cal M}_{SBS}(c_1(L))$ and simultaneously illustrate
 the correspondence "special Bohr - Sommerfeld cycles = Hamiltonian isotopy classes of D - exact lagrangian submanifolds".

 \section{Geometrical interpretations}

 Consider  compact simply connected symplectic manifold $(M, \omega)$ of real dimension $2n$. Suppose that the cohomology class $[\omega]$ is integer.
 Then it exists a complex line bundle $L \to M$ such that $c_1(L) = [\omega]$, called prequantization bundle, so fixing a hermitian structure on it
 one gets the space of hermitian connections ${\cal A}_h(L)$ with a distinguished subset ${\cal O}(\omega)$ consists of such $a$ that $F_a = 2 \pi \omega$.
 Since our $M$ is simply connected this subset ${\cal O}(\omega)$ is an orbit of the gauge group action. For any smooth lagrangian submanifold $S \subset M$
 and any connection $a \in {\cal O}(\omega)$ the restriction $(L, a)|_S$ is a flat line bundle, and we say that $S$ is Bohr - Sommerfeld iff
 $(L, a)|_S$ admits a covariant constant section $\sigma_S$ uniquely defined up to scaling. The BS - property on $S$ does not depend on the particular choice
 of  $a \in {\cal O}(\omega)$ while $\sigma_S$ does. Indeed, for any other $a_1 \in {\cal O}(\omega)$ the difference $\nabla_{a_1} - \nabla_a =
 \imath d \phi$ where $\phi \in C^{\infty}(M, \mathbb{R})$, and new covariantly constant section $\sigma^1_S$ reads as $e^{-\imath \phi} \sigma_S$.

 Thus the definition of the subset ${\cal U}_{SBS} \subset \mathbb{P} \Gamma (M, L) \times {\cal B}_S$, where ${\cal B}_S$ is the moduli space of Bohr - Sommerfeld
 lagrangian submanifolds, depends on the choice of the prequantization
 connection $a \in {\cal O}(\omega) \subset {\cal A}_h (L)$ as well  since our speciality condition
 $$
 ([\alpha], S) \in {\cal U}_{SBS}(a) \quad \Longleftrightarrow \quad \frac{\alpha|_S}{\sigma_S} = e^{\imath c} f, \quad f \in C^{\infty}(S, \mathbb{R}_+),
 $$
   depends on the covariantly constant section $\sigma_S$ of the restriction $(L, a)|_S$.

  For two different connections $a_2, a_1 \in {\cal O}(\omega)$ the corresponding subsets ${\cal U}_{SBS} (a_i)$ intersect each other at the following subset:
  $$
  {\cal U}_{SBS}(a_2) \cap {\cal U}(a_1) = \{ (p, S) \in {\cal U}_{SBS}(a_1) \quad \vert \quad \phi|_S = \rm{const} \},
  \eqno 1
  $$
  where $\phi \in C^{\infty}(M, \mathbb{R})$ is given by $\nabla_{a_2} - \nabla_{a_1} = \imath \rm{d} \phi$. Indeed, since by the very definition
  $\alpha \in \Gamma(M, L)$, corresponding to $p = [\alpha] \in \mathbb{P} \Gamma(M, L)$, never vanishes on $S$, the pair $(p, S)$
  can belong to both ${\cal U}_{SBS} (a_i)$ if and only if $\sigma_S^1 = C \sigma^2_S$ therefore $\rm{d} \phi|_S$ must be trivial.

  Recall that each smooth function $F$ on $M$ generates certain vector field $\Theta(f)$ on ${\cal B}_S$: at point $S$ its value is given by the exact form
  $\rm{d} (F|_S)$, and the flow, generated by $\Theta(F)$, is precisely given by the flow  generated by the Hamiltonian vector field $X_F$,
   applied to lagrangian submanifolds (the details can be found in [5], [7]). Therefore  the intersection (1) can be geometrically described as follows:
  the difference function $\phi$ generates the corresponding vector field $\Theta(\phi) \in \rm{Vect} {\cal B}_S$, and if we take the second canonical projection
  $  p_2: {\cal U}_{SBS}(a_1) \to {\cal B}_S$ then the intersection (1) is given by the preimage $p_2^{-1}((\Theta(\phi))_0)$ of the zeroset of
  $\Theta(\phi)$.

  Since the intersection condition is imposed on the second element of the pair $(p, S)$ only it hints some geometrical description of ${\cal U}_{SBS}(a)$.

  Both the direct summands $\mathbb{P} \Gamma(M, L)$ and ${\cal B}_S$ carry natural $U(1)$ - bundles. For the first direct summand it is presented
  by ${\cal O}(1)$, the standard line bundle over the projective space. Since the Kahler structure on the projective space is fixed, this bundle
  carries the corresponding hermitian connection $A$ with the curvature form $F_A = 2 \pi \Omega_{FS}$. Note that the original $U(1)$ action on the prequantization
  bundle $L \to M$ generates the corresponding action on ${\cal O}(1)$.

  On the other hand the moduli space ${\cal B}_S$ carries a natural $U(1)$ - bundle ${\cal P}_S(a) \to {\cal B}_S$ (see [5]) formed by the Planckian cycles
  in the contact manifold $\rm{tot} (S^1(L) \to M)$. The fiber over $S \in {\cal B}_S$ consists of covariantly constant lifting of $S$ to
  $(S^1(L),a)|_S$, and the $U(1)$ action is induced again by the original $U(1)$ action on $L \to M$; note that the bundle ${\cal P}_S(a) \to {\cal B}_S$
  depends on the choice of $a \in {\cal O}(\omega)$ since its fibers do.

  Thus on the direct product $\mathbb{P} \Gamma (M, L) \times {\cal B}_S$ one has two $U(1)$ - bundles $p_1^* {\cal O}(1)$ and $p_2^* {\cal P}_S(a)$.
  Then the subset ${\cal U}_{SBS}(a)$ in the direct product can be characterized by the following fact:

  {\bf Proprosition 1.}  {\it The bundles $p_1^* {\cal O}(1)$ and $p_2^* {\cal P}_S(a)$ are canonically isomorphic to each other being restricted to ${\cal U}_{SBS}(a)$.}

  Indeed, the fiber of the first bundle over a point $(p, S)$ is given by elements $e^{it} \alpha^*$ pointwise dual to the section $\alpha$ such that $p = [\alpha]$.
  At the same time the fiber of the second bundle is given by $e^{it} \sigma_S$; therefore in the direct product of the fibers $U(1) \times U(1)$ one has
  some diagonal $U(1)$ distinguished by the condition $\alpha^*|_S(\sigma_S) \in C^{\infty}(S, \mathbb{R}_+)$ --- the natural paring is real strictly positive
  along $S$.  It happens exactly over ${\cal U}_{SBS}(a)$ due to the speciality condition.

  Note that we have even more: the subset ${\cal U}_{SBS}(a)$ is the biggest possible subset where the lifted bundles are canonically isomorphic to each other.

  {\bf Corollary.}  {\it The subset ${\cal U}_{SBS}(a)$ carries a natural $U(1)$ - bundle ${\cal L}$ which is isomorphic to the restrictions of
  $p_1^*{\cal O}(1)$ and $p_2^*{\cal P}_S(a)$.}

  In particular bundle ${\cal L} \to {\cal U}_{SBS}(a)$ carries hermitian connection $\tilde A = p_1^* A$ with the curvature form $F_{\tilde A} = p_1^* \Omega_{FS}$
  given by the weak Kahler structure form. Therefore any other hermitian connection on ${\cal L}$ is given by the corresponding 1 - form on ${\cal U}_{SBS}(a)$,
  and a natural question here is to correct $\tilde A$ by an appropriate 1 - form, lifted from ${\cal B}_S$, such that the curvature form of the resulting connection
  would be strong Kahler or symplectic form.

  \section{Transformations}

  At the same time it is not hard to see that all ${\cal U}_{SBS}(a)$ are isomorphic to each other. Indeed, for any smooth real function $F \in C^{\infty}(M, \mathbb{R})$,
  globally defined over $M$, one has the corresponding projective transformation $P(F) \in \rm{Aut} \mathbb{P} \Gamma (M, L)$ given by formula $[\alpha] \mapsto
  [e^{i F} \alpha]$. Since on the level of vector space $\Gamma(M, L)$ this transformation preserves the hermitian scalar product
  $$
  \int_M <e^{iF} \alpha_1, e^{i F} \alpha_2>_h d \mu_L = \int_M <\alpha_1, \alpha_2>_h d \mu_L,
  $$
  the transformation $P(F)$ is a Kahler isometry of $(\mathbb{P}\Gamma (M, L), \Omega_{FS})$.

   For two hermitian connections $a_1, a_2$ such that the difference $\nabla_{a_2} - \nabla_{a_1} = \imath \rm{d} F$ for a smooth function $F \in C^{\infty}(M, \mathbb{R})$
   one has that $([\alpha], S) \in {\cal U}_{SBS}(a_1)$ if and only if $(P(F)([\alpha]), S) \in {\cal U}_{SBS}(a_2)$. Indeed,  $\alpha|_S = e^{i c} f \sigma_S^1
   = e^{i c} f e^{-i F|_S} \sigma_S^2$ if and only if $e^{i F} \alpha|_S = e^{i c} f \sigma_S^2$, therefore:

   {\bf Proposition 2.} {\it  All the subspaces ${\cal U}_{SBS}(a)$ are isomorphic to each other; the isomorphism is given by  an appropriate
   transformation $P(F)$.}

   Here we strongly exploit the fact that the first components belong just to a projective space. More complicated question is about the variation of
   the second components: Bohr - Sommerfeld lagrangian submanifolds.

  The set $\mathbb{F} = \{ P(F), F \in C^{\infty}(M, \mathbb{R} \}$ is an abelian subgroup in $\rm{Aut} \mathbb{P} \Gamma (M, L)$, and
  it is not hard to see that there exist the projective subspaces which are invariant with respect to this subgroup. Let $D \subset M$ be a submanifold which represents class
  $P.D. [\omega] \in H_{2n-2}(M, \mathbb{Z})$ then the  subspace $\mathbb{P}(D) = \{ [\alpha] \vert (\alpha)_0 = D \} \subset \mathbb{P} \Gamma (M, L)$ is exactly of this type.
 Note that $D$ can be non smooth, consisting of several components with multiplicities. The subspace $\mathbb{P}(D)$ is not projective since we require in the definition
 that the zeroset of $\alpha$ is exactly $D$; however as we will show below it is an affine space.

 We call the  group $\mathbb{F}$ {\it phase changing group} since its elements change the phases of sections.

  The phase changing group contains one - parameter subgroups $\{ P(tF) \}$, therefore each smooth function $F$ generates vector field $\Theta_P(F)$ as the infinitesimal part of
  $P(tF)$, and by the very definition this vector field preserves the Kahler structure on $\mathbb{P} \Gamma (M, L)$. Since $\mathbb{F}$ is commutative, the space of all
  such $\Theta_P(F)$ gives us an integrable distribution on the projective space $\mathbb{P} \Gamma (M, L)$. It is not hard to see that this distribution is integrable.
  Since for a section $\alpha \in \Gamma (M, L)$ the transformation $\alpha \mapsto e^{i F} \alpha$ preserves the pointwise norm $\vert \alpha \vert_h$ then
   the projective version $P(F)$ must preserve real 1 -form $\rm{d} \rm{ln} \vert \alpha \vert_h$, correctly defined on the complement $M \backslash D$.

    Recall from [4] that for the component $\mathbb{P}(D)$ one has the following attachment: every class $[\alpha] \in \mathbb{P}(D)$ is presented by the corresponding
  complex 1 - form $\rho(\alpha) = \frac{\nabla_a \alpha}{\alpha}$ correctly defined on the complement $M \backslash D$, such that $\rm{Re}(\rho(\alpha))$ is exact and
  $\rm{d} \rm{Im} (\rho(\alpha)) = 2 \pi \omega$ (the real part is just $\rm{d} \rm{ln} \vert \alpha \vert_h$). The point is that it is one - to one correspondence:
  every complex 1- form $\rho$ of this type gives a section $\alpha$ vanishing along $D$ defined uniquely up to scaling, see [4]. In this presentation the action of $P(F)$
  looks  quite simple: $\rho(P(F)(\alpha)) = \rho(\alpha) + \imath \rm{d} F$. And the crucial fact has been established in [1]: in this presentation SBS condition
  reads just as $\rm{Im} \rho(\alpha)|_S \equiv 0$.

   Therefore the  subspace $\mathbb{P}(D)$ is an affine space associated with the complex vector space formed by exact complex 1 -forms on the complement $M \backslash D$.
    This affine space is fibered by the attachment $[\alpha] \mapsto \rm{Re} \rho(\alpha)$, and
    if we denote as $\mathbb{P}^0 (D)$ the fiber over a given point $\rm{Re} \rho(\alpha_0)$ then the vector field $\Theta_P(F)$ must be tangent to
    this subspace, and one has

    {\bf Proposition 3.} {\it The distribution spanned by vector fields $\Theta_P(F)$ is integrable: its leaves are given by subspaces of the form $\mathbb{P}^0 (D)$.}

    Note that $\mathbb{P}^0(D)$ is again an affine space but it is real one, associated with real vector space of exact 1 - forms. But for our SBS story it is not important since

    {\bf Proposition 4.} {\it If for a pair $(p, S) \in {\cal U}_{SBS}(a)$ the first element $p \in \mathbb{P}(D)$ then it exists some pair $(p', S)
    \in {\cal U}_{SBS}(a)$ such that $p' \in \mathbb{P}^0 (D)$.}

    Indeed, SBS condition does not depend on the real part $\rm{Re} \rho(\alpha)$ therefore we can reduce $p$ to $p'$ just adding  appropriate real exact form
    to $\rho(\alpha)$.

    Therefore to study ${\cal U}_{SBS}(a)$ one needs just to study the situation over components $\mathbb{P}^0(D)$ not over whole $\mathbb{P} \Gamma (M, L)$. In particular
    it is interesting to find certain universal rule how to choose the real part $\rm{Im} \rho$ for a given $D \subset M$. In some cases the choice can be canonically made:
    suppose that our symplectic manifold $(M, \omega)$ admits an integrable complex structure $I$ which is compatible with $\omega$. Thus our $(M, \omega, I)$
     is a complex Kahler variety with the Kahler structure of the Hodge type. The choice of a hermitian connection $a_I$ in the orbit ${\cal O}(\omega)$
    induces a holomorphic vector bundle structure on $L$, and it is well known that the corresponding space of holomorphic sections $H^0(M_I, L)$ is a finite dimensional
    subspace in $\Gamma (M, L)$. Every holomorphic section $\alpha \in H^0(M_I, L)$ modulo scaling corresponds to its zero divisor $D_{\alpha}$ which is formed by
    complex submanifolds with multiplicities. Thus the projectivized space $\mathbb{P} H^0(M_I, L)$ is called the complete linear system $\vert L \vert$.

    For this case the attachment
    $$
    D_{\alpha} \leftrightarrow [\alpha] \leftrightarrow \rm{Re} \rho(\alpha) \leftrightarrow \mathbb{P}(D_{\alpha})
    $$
    is correct, and we have the following fact: the intersection $\mathbb{P} H^0(M_I, L) \cap \mathbb{P}(D) \subset \mathbb{P} \Gamma (M, L)$ is at most
    one single point, and it is non trivial if and only if $D$ is formed by complex submanifolds with multiplicities so $D \in \vert L \vert$.

    Thus for this case we have a finite dimensional set of subspaces $\{ \mathbb{P}(D)^0 \quad \vert \quad D \in \vert L \vert \}$ with marked points $p_D \in
    \mathbb{P}(D)^0$. Note that these $\mathbb{P}^0(D)$ do not intersect each other by the very definition.

       \section{Deformations}

   In the previous section we study possible transformations of the first elements of our pairs $(p, S)$; now we have to study
   the deformations of the second ones.

  {\bf Proposition 5.}  {\it Let $([\alpha], S)$ be a  point in ${\cal U}_{SBS}(a)$ for a fixed $a \in {\cal O}(\omega)$. Then  for each
   small Bohr - Sommerfeld lagrangian variation $S_{\delta}$ of given $S$ it exists
  the corresponding deformation pair $(P(F_{\delta}[\alpha], S_{\delta})$  which is again  contained by ${\cal U}_{SBS}(a))$.}

  Indeed, if $S$ is special Bohr - Sommerfeld with respect to $\alpha$ then for a Darboux - Weinstein neighborhood ${\cal O}_{DW}(S) \subset M$
  of $S$ one has the corresponding 1-  form $\frac{1}{2\pi} \rm{Im} \rho(\alpha)$ such that its differential is $\omega$ and its restriction to $S$
  identically vanishes. Take the canonical 1- form $\alpha_{can}$ on ${\cal O}_{DW}(S)$ and consider the difference form $\alpha_{can} - \frac{1}{2\pi} \rm{Im} \rho(\alpha)$.
  It is closed but since the neighborhood ${\cal O}_{DW}(S)$ can be contracted to $S$ where this form identically vanishes by the assumptions it exists a smooth function
  $F_0$ such that $\rm{d} F_0 = \alpha_{can} - \frac{1}{2\pi} \rm{Im} \rho(\alpha)$ over ${\cal O}_{DW}(S)$. Note that at the same time by the definition $\rm{d} F_0|_S \equiv 0$.

  Consider new hermitian connection $a_0$ such that $\nabla_{a_0} = \nabla_a - \imath \rm{d} F_0$; it is clear that $([\alpha], S)$ still belongs to ${\cal U}_{SBS}(a_0)$.
  At the same time the local picture near $([\alpha], S)$ in the last space looks as follows. Since there $\alpha_{can} \equiv \frac{1}{2 \pi} \rm{Im} \rho(\alpha)$
  then for any small deformation $S_{\delta}$ of $S$, given by a smooth function $\phi \in C^{\infty}(S, \mathbb{R})$ as usual for the Darboux - Weinstein
  presentation (see [GT]), one has $\frac{1}{2 \pi} \rm{Im} \rho(\alpha)|_{S_\delta} = \rm{d} \pi^* \phi$ (here $\pi: {\cal O}_{DW}(S) \to S$ is the canonical projection).
  Therefore if one takes a variation $S_{\delta}$ of our given $S$, represented in ${\cal O}_{DW}(S)$ by the corresponding smooth function $\phi \in C^{\infty}(S, \mathbb{R})$
  the restriction $\frac{1}{2\pi} \rm{Im} \rho(\alpha)|_{S_{\delta}}$ equals $\rm d(\pi^* \phi)|_{S_{\delta}}$ where $\pi: {\cal O}_{DW}(S)$ is generated by the canonical projection $T^*S \to S$. The real function $\pi^* \phi|_{S_{\delta}}$ can be extended to  ${\cal O}_{DW}(S)$ and then to whole $M$, and we denote this extension as $F_1 \in C^{\infty}(M, \mathbb{R}$.

  Now if one takes the space ${\cal U}_{SBS}(a_0 + \imath \rm{d} F_1)$ then  the pair  $([\alpha], S_{\delta})$ evidently satisfies SBS condition
  with respect to this deformed hermitian connection. Consequently the same pair $([\alpha], S_{\delta})$ belongs to ${\cal U}_{SBS}(a+ \imath \rm{d}(F_0+ F_1))$,
  and if we put $\delta = F_0 + F_1$ then $(P(\delta) [\alpha], S_{\delta})$ must seat inside ${\cal U}_{SBS}(a)$, and it ends the proof.

  Note that during the construction we did not leave the same projective subspace $\mathbb{P}(D)$.

  These arguments lead to the following:

  {\bf ${\cal B}_S$- Covering theorem.} {\it Let $S_t, t \in [0;1]$, be a Hamiltonian isotopy of a Bohr - Sommerfeld lagrangian submanifold $S_0$ such that $([\alpha_0], S_0)
  \in {\cal U}_{SBS}(a)$ for a fixed connection $a$. Suppose that for all $t \in [0;1]$  the intersection $S_t \cap D_{\alpha_0} = \emptyset$ and that $S_t$ is smooth. Then
  it exists the corresponding family $([\alpha_t], S_t) \in {\cal U}_{SBS}(a)$ such that all $[\alpha_t]$ belongs to the same $\mathbb{P}^0(D_{\alpha_0})$.}

  Indeed, since the Hamiltonian isotopies preserve the Bohr - Sommerfeld condition, every $S_t$ is BS; and for each $t \in (0;1)$ the corresponding
  $S_t$ has a Darboux - Weinstein neighborhood ${\cal O}_{DW}(S_t)$ which contains every $S_{t'}$ for $t' \in (t- \varepsilon, t+ \varepsilon)$.
  The segment $[0;1]$ is compact therefore it exists a finite choice of
  $S_{t_i}, i = 1, ... , N$, such that the union of  neighborhoods $\bigcup_{i=1}^N {\cal O}_{DW}(S_{t_i})$ contain every $S_t$. Then we can apply
  Proposition 5 several times, passing through the segment, and establish that $S_1$ carries the corresponding class $[\alpha_1]$ such that
  the statement of ${\cal B}_S$-  Covering theorem above holds.

  At the same time we can use another type of arguments, based on   the stability of SBS property with respect to flows generated by Hamiltonian vector fields.
   Namely let $F$ be a global function on $M$, $X_F$ is its Hamiltonian vector
  field and $\phi^t_{X_F}$ is the flow generated by $X_F$. Then for a pair $(p, S) \in {\cal U}_{SBS}(a)$ one has the corresponding Hamiltonian deformation
  $(\phi^t_{X_F}(p), \phi^t_{X_F} (S))$ and the point is that

  {\bf Proposition 6.} {\it The pair  $(\phi^t_{X_F}(p), \phi^t_{X_F} (S))$ belongs to ${\cal U}_{SBS}(a)$.}

  We prove this statement in the following version: let the zero set  $D$ for the class $p$ be stable with respect to the flow $\phi^t_{X_F}$. Then
  the action of $\phi^t_{X_F}$ on the pair $(p, S)$ can be reformulated as the action of the flow on the pair $(\rho, S)$ where complex 1 -form
  $\rho$ corresponds to $p$ on $M \backslash D$. Since the attachment $p \leftrightarrow \rho$ is correct under the condition that $a$ and $D$ are fixed,
  we can forget about sections and work with 1 - forms. Then the flow $\phi^t_{X_F}$ gives as a family $(\rho_t, S_t)$, and obviously the condition
  $\rm{Im} \rho|_S = 0$ is stable with respect to the deformations.

\section{Definition of the moduli space}

Consider a  simply connected smooth compact (or projective) algebraic variety $X$ together with a very ample line bundle $L \to X$ which does exists by the very definition.
Choose an appropriate hermitian structure $h$ on $L$ and take the corresponding Kahler form $\omega_h$ given by Kahler potentials $\psi_{\alpha} =
- \rm{ln} \vert \alpha \vert_h$ for holomorphic sections $\alpha \in H^0(X, L)$ on the complements $X \backslash D_{\alpha}$. At the same time
the choice of $h$ in presence of the holomorphic structure on $L$ gives a distinguished connection $a_I \in {\cal O}(\omega_h)$.

Thus the choice of $h$ leads to the situation which we have studied above: one gets $(M = X, \omega = \omega_h, I, L, a_I)$, and we can apply the constructions
of SBS geometry. At the same time we know, see [2], that for any $[\alpha] \in \mathbb{P} H^0(X, L) \subset \mathbb{P} \Gamma (X, L)$ the preimage
$p_1^{-1}([\alpha] \in {\cal U}_{SBS}(a_I)$ consists of the following terms. Take the corresponding Kahler potential $\psi_{\alpha}$ on
$X \backslash D_{\alpha}$, take its critical points $x_1, ..., x_N$ and the corresponding finite trajectories of the gradient vector field
$\rm{grad} \psi_{\alpha}$, joining $x_i$'s; then the union of the finite trajectories gives the Weinstein skeleton $W(X \backslash D_{\alpha})$
of the complement, and a smooth BS lagrangian submanifold $S \subset X \backslash D_{\alpha}$ is special with respect to $[\alpha]$ if and only if
$S \subset W(X \backslash D_{\alpha})$ (the details can be found in [2]). On the other hand, see [6], even in the simplest cases $W(X \backslash D_{\alpha})$
does not admit smooth components therefore $p^{-1}_1([\alpha])$ is empty set.

{\bf Example ([3]):} consider $X = \mathbb{C} \mathbb{P}^1$, $L = {\cal O}(3)$, then for generic section the Weinstein skeleton $W(X \backslash D_{\alpha})$
is presented by three valent graph on the 2 - sphere with three vertices and three edges. Therefore no smooth closed loops exist there, and the preimage
is empty set.

On the other hand this emptiness follows from the fact that the situation is very far to be generic: the chosen points $[\alpha]$ and the defining connection
$a_I$ are strongly related to each other. If we use a small perturbation $\delta: \mathbb{P} H^0(X, L) \to C^{\infty}(X \backslash D_{\alpha}, \mathbb{R})$
such that $[\alpha] \mapsto P(\delta([\alpha])) [\alpha]$ in the same affine space $\mathbb{P}(D_{\alpha})$ one can expect that the corrected preimage
$p_1^{-1}(P(\delta ([\alpha])))$ should be non empty. The term "small" can be understood as follows: since each function
 $\delta([\alpha])$ is smooth globally defined on compact manifold $X$, we can use the universal bound $\rm{max} \delta(p) - \rm{min} \delta(p) \leq \varepsilon$
 for each $p$.  If we denote as $\mathbb{P}H^0(X, L)_{\delta} \subset \mathbb{P} \Gamma (X, L)$ the corresponding deformation
of the projective space, then it is not hard to see that this deformation gives

{\bf Proposition 7.} {\it For a sufficiently  small deformation $\delta: \mathbb{P} H^0(X, L) \to C^{\infty}(X \backslash D_{\alpha}, \mathbb{R})$
the resulting space $\mathbb{P} H^0(X, L)_{\delta}$ is a smooth real $2 (h^0(X, L) -1)$  submanifold, symplectic with respect to $\Omega_{FS}$.}

  Indeed, as we have seen above the projective space $\mathbb{P} H^0(X, L)$ is associated to the corresponding family $\{ \mathbb{P}^0 (D_{\alpha}) \}$
  of  affine spaces which do not intersect each other and which are transversal to the projective space. The transformations generated by $\delta$
  act along "slices" $\mathbb{P}(D_{\alpha})$, and since the deformation is sufficiently small the resulting submanifold $\mathbb{P} H_0(X, L)_{\delta}$
  must be symplectic.

  Now let $(X, L)$ be as above. For an appropriate hermitian structure $h$ we define

  {\bf Definition.}{\it The moduli space of special Bohr - Sommerfeld cycles
  $${\cal M}_{SBS}(c_1(L), \rm{top} S, [S]) = p_1^{-1}(\mathbb{P}H^0(X, L)_{\delta}) \subset {\cal U}_{SBS}(a_I)
  $$
   where $\delta$ is a generic sufficiently small deformation, $\rm{top} S$ --- topological type of S and $[S] \in H_n(X, \mathbb{Z})$ a fixed homology class.}

   The dependence on $h$ has been discussed in [2]; our aim now is to show that the geometry of ${\cal M}_{SBS}$ does not depend on the choice of the small deformation.

   {\bf Proposition 8.} {\it The space ${\cal M}_{SBS}$ does not depend on the choice of generic deformation.}

   We can check this fact separately for the deformations of a given point $[\alpha] \in \mathbb{P}H^0(X, L)$; the corresponding preimage is non trivial "in general"
   if for any generic deformation $\delta([\alpha])$ it is non trivial therefore we must compare two sets $\{S_1^1, ...., S_l^1, ... \}$ and $\{S_1^2, ..., S_k^2, ... \}$
   such that $(\delta_i([\alpha]), S^i_j) \in p^{-1}_1(\mathbb{P} H_0(X, L)_{\delta_i}$ where $i = 1, 2$ taking into account that the deformations $\delta_i$ are sufficiently small.

   Note that on the complement $X \backslash D_{\alpha}$ we have three vector fields: $\lambda = \omega^{-1}(\rm{Im} \rho(\alpha)), \lambda_i = \lambda + X_{\delta_i}$,
   $i = 1, 2$, where $X_{\delta_i}$ denotes the Hamiltonian vector field of the function $\delta_i([\alpha])$. It is clear that $\lambda_i$ is exactly $\omega^{-1}( \delta_i([\alpha])$,
   therefor the SBS - condition for the vector fields reads as follows: $\lambda_i$ is tangent to smooth Bohr - Sommerfeld lagrangian submanifold $S_j^i$ for each $j$. Indeed
   if any 1 -form $\rho$ satisfies $\rho|_S \equiv 0$ for a lagrangian submanifold $S$ then the vector field $\omega^{-1}(\rho)$ must be tangent to $S$.

  In general the last fact does not imply that $\omega^{-1}(\rho)|_S$ admits too many zeros, but in our situation we have the following remark: since $\lambda$
  is the gradient vector field of the Kahler potential $\psi_{\delta}$ it admits zeros coming with the critical points of a function (so the number is dictated by
  the Morse inequality); at the same time our deformations $\delta_i$ are sufficiently small therefore we can expect that the vector fields $\lambda_i$ admit
  the same numbers of zeros, and by the observation made above these zeros lie on $S^i_j$.

  Now we can join the vector fields $\lambda_1$ and $\lambda_2$ by a path $\{\lambda_t | t \in [1;2] \}$ of vector fields which satisfy the same property about the zeros.
  Therefore it induces certain correspondence between components $S^1_j$ and $S^2_l$: namely take the zeros of $\lambda_1$ on the first component $S_1^1$
  and then follow how the corresponding zero points move along the path $\{ \lambda_t \}$ when $t$ goes from 1 to 2. Using this attachment we establish
  that the number of components $S^1_j$ must be the same as for $S^2_i$. Globalizing further over the complete linear system we get an identification of $p^{-1}_1(\mathbb{P} H^0(X, L)_{\delta_1})$  and $p^{-1}_1 (\mathbb{P} H^0(X, L)_{\delta_2})$, and it ends the proof.

  {\bf Remark.} Of course, here we exploit certain very special version of a statement which can be called {\bf $\mathbb{P}$ - covering theorem}. Above
  we present and prove ${\cal B}_S$ - covering theorem when a deformation of points on the moduli space ${\cal B}_S$ is covered by the corresponding deformation
  in ${\cal U}_{SBS}(a)$. Now if we take a deformation of other elements which live on the projective space $\mathbb{P} \Gamma (X, L)$ then the story in a sense
  turns to be even simpler: since the space $\mathbb{P}(D_{\alpha})$ is affine then every two elements $p_1, p_2 \in \mathbb{P} (D_{\alpha})$ can be joined just by
  the corresponding segment. Indeed, having the  distinguished point $[\alpha]$ each other element $p_i \in \mathbb{P}(D_{\alpha})$ is uniquely represented by
  1 - form $\rm{Im} \rho(\alpha) + \rm{d} F$ where $F$ is a smooth function on $X$, therefore the segment between $p_1$ and $p_2$ is presented by
  the family $\rm{Im} \rho(\alpha) + t \rm{d} F_1 + (1-t) \rm{d} F_2$; moreover, if this linear deformation is not suitable for a problem then
  we can take any path in the space $\rm{d} \Omega^0_X$  of exact 1- forms with ends at $\rm{d} F_1$ and $\rm{d} F_2$. But evidently such $\mathbb{P}$ - covering
  theorem can not be true for every path as it was for ${\cal B}_S$ - covering theorem: for example, take $F_2 = -F_1$, then for the middle point $p
  = [\alpha]$ no smooth SBS lagrangian submanifolds appear as we have mentioned above, and at the same time for the end points the corresponding SBS - submanifolds
  can exist. We will discuss  details related to this question in the last section.

  \section{Exact lagrangian submanifolds}

  In [3] for absolutely the same situation --- algebraic variety $X$ and very ample line bundle $L$ --- one constructs certain moduli space of D - exact lagrangian
  submanifolds, denoting this space as $\tilde {\cal M}_{SBS}$, and then derives in this space a stable component ${\cal M}_{st} \subset \tilde {\cal M}_{SBS}$, promising explanations why
  the first moduli space is labeled by the same letters "SBS" and how stable D - exact lagrangian submanifolds are related to SBS lagrangian submanifolds.

  Briefly, for a simply connected projective algebraic variety $X$ and a very ample line bundle $L \to X$ one fixes an appropriate hermitian structure $h$, endowing
  $X$ by the corresponding symplectic form $\omega_h$. Then to a point $p = [\alpha] \in \mathbb{P} H^0(X, L)$ one attaches the set of D - exact smooth lagrangian submanifolds
  $\{ S \vert S \subset X \backslash D_{\alpha} \}$ of the same topological type $\rm{top} S$ and the same homology class $[S] \in H_n(X, \mathbb{Z})$ such that
  $S$ represents a non trivial homology class in $H_n(X \backslash D_{\alpha}, \mathbb{Z})$; then the quotient space by the Hamiltonian isotopies on $X \backslash D_{\alpha}$
  gives a discrete set, therefore globalizing over whole $\mathbb{P} H^0(X, L)$ one gets a space $\tilde {\cal M}_{SBS}$. In [3] one proves that this space
  is an open Kahler manifold, fibered over an open subset of $\mathbb{P} H^0(X, L)$ without ramifications. The notion of D - exactness essentially coincides with the
  standard notion of exactness: a lagrangian submanifold is called to be exact if the restriction of $\rm{Im} \rho(\alpha)|_S$ is exact; evidentally then the same happens
  for each element from $\mathbb{P}(D_{\alpha})$. One uses D - exactness in [3] for the following reason: in this formulation the stability of exactness with respect
  to deformations of $D_{\alpha}$ is much more clear: we can vary $[\alpha] \in \mathbb{P} H^0(X, L)$ with a fixed exact $S$ and before $D_{\alpha}$ non trivially  intersects $S$
  the last one has to be exact for any deformed $[\alpha]$. Below we do not study deformations of this type and therefore the standard definition of exactness is exploited.

  In [3] one outlines a connection between "exact" construction and SBS - construction based on the following suggestion: we say that a cycle $\Delta $ sitting inside
  the Weinstein skeleton $W(X \backslash D_{\alpha}$ admits a Bohr - Sommerfeld resolution if it exists a homotopy $\{ S_t \}, t \in [0;1],$ such that
  $S_0 = \Delta$, and for each $t \in (0;1]$ the corresponding $S_t \subset X \backslash D_{\alpha}$ is a smooth Bohr - Sommerfeld lagrangian submanifold.
  From the definition it follows that $S_t$ represents a non trivial class in $H_n(X \backslash D_{\alpha}, \mathbb{Z})$; moreover from the stability of
  exactness it follows that each $S_t$ must be exact. Thus we can derive from the moduli space $\tilde {\cal M}_{SBS}$ of exact lagrangian submanifolds
  the components ${\cal M}^{st}$ formed by the classes, which admit representatives given by Bohr - Sommerfeld resolutions.

  Now one has the following

  {\bf Theorem.} {\it The moduli space ${\cal M}_{SBS}$ defined above is naturally isomorphic to the component ${\cal M}^{st}$ of the moduli space of D - exact
  lagrangian submanifolds.}

  The geometrical essence of this statement is the following: if we have smooth SBS lagrangian submanifold $S_{\delta}$ for generic small deformation $\delta$ then
  this $S_{\delta}$ must be exact since it is Bohr - Sommerfeld and close to the Weinstein skeleton $W(X \backslash D_{\alpha})$ (we can think that outside of a certain
  neighborhood of $W(X \backslash D_{\alpha}$ the correction term $\rm{d} \delta$ is trivial). Then taking general deformation family $\delta_t$ where $t$ tends to zero
  we can find the desired homotopy. On the other hand, if an appropriate homotopy $\{ S_t \}$ does exist then the restriction $\rm{Im} \rho (\alpha)|_{S_t}$ is presented by
  an exact 1 - form $\rm{d} \phi_t$ on $S_t$, and we can extend the corresponding function $\phi_t$ to a small neighborhood of $S_t$ such that outside of this neighborhood
   this extension is constant. This construction gives us  the desired deformation $\delta_t$ such that $\delta_t|_{S_t} = \phi_t$;  it is clear that this deformation is small.
   This ends the proof.

   {\bf Remark.} At this point we would like to discuss a natural question: why we need this SBS - constructions if we already have certain finite dimensional moduli space presented by
   D - exact lagrangian submanifolds? The main reason is the following: in Section 1 above we present the construction of universal bundle ${\cal L}  \to {\cal U}_{SBS}(a)$,
   by the very definition, see Section 4, the moduli space ${\cal M}_{SBS}$ is embedded to ${\cal U}_{SBS}(a_I)$ therefore it exists the restriction ${\cal L}_{\delta}
   \to {\cal M}_{SBS}$ which depends in principle  on the choice of the small deformation $\delta$. However it is not hard to establish that topologically
   this restriction ${\cal L}_{\delta}$ does not depend on the choice of $\delta$, therefore we get much more interesting object than just a manifold ${\cal M}_{SBS}$ ---
   a pair "mainfold + bundle".  The variation of $\delta$ on the level of the bundle corresponds to the choice of  hermitian structure: the $\delta$ - variation can be reformulated
   in terms of the deformation of the basic connection $a_I$, and the last one as we have seen above corresponds to the changing of basic covariantly constant sections $\sigma_S$,
   and the last one corresponds to the changing of normalized frame at a given point which is of course reflected by the corresponding changing of hermitian structure.
   On the other hand in [3] one states a conjecture: the component ${\cal M}^{st}$ of the moduli space of D - exact lagrangian submanifolds is isomorphic to the complement
   of an algebraic variety to an ample divisor. This conjecture can be studied in terms of the pair ${\cal L} \to {\cal M}_{SBS}$: first, find a Kahler structure on
   the base; second, find a hermitian connection on ${\cal L}$ such that its curvature form is proportional to the Kahler form; third, find an appropriate section
   such that its zeros does not belong to ${\cal M}_{SBS}$. Realization of this programme requires many technical details, starting with the question about
   connections on ${\cal L}$ posted in Section 1; even if the conjecture is not true the investigations in differential geometry of the moduli space ${\cal B}_S$ of Bohr - Sommerfeld
   lagrangian submanifolds shall be useful for possible application in Geometric Quantization.

   \section{Weinstein structures and Eliashberg conjectures}

   Special Bohr - Sommerfeld geometry is closely related to the theory of Weinstein structures. Recall, see [6], that a vector field $\lambda$ on an open symplectic  manifold $M \backslash D_{\alpha}$ with symplectic form $\omega$ is called Liouville if the Lie derivative ${\cal L}_{\lambda} \omega = \omega$. As we have seen for any regular section
   $\alpha \in \Gamma (M, L)$ with zeroset $D_{\alpha}$, presented by a combination of smooth $2n - 2$ components with multiplicities, one has the corresponding vector field
   $\lambda_{\alpha} = \omega^{-1}(\rm{Im} \rho(\alpha)$ which is Liouville. By the definition this vector field depends on the class $[\alpha] \in \mathbb{P}\Gamma (M, L)$ only. At the same time it is defined a real function $\psi_{\alpha} = - \rm{ln} \vert \alpha \vert_h$ with pole along $D_{\alpha}$. Note however that the space of possible Lioville fields
   on $M \backslash D_{\alpha}$ is not exhausted by the attachment: since $H^1(M \backslash D_{\alpha}, \mathbb{Z})$ is non trivial (it follows from the fact that the prequantization
   bundle is topologically nontrivial), then $H^1(M \backslash D_{\alpha}, \mathbb{Z})$ must be non trivial too, therefore one can add a closed but  non exact 1 - form to
   $\rm{Im} \rho(\alpha)$ and then apply $\omega^{-1}$ to the sum. Therefore SBS - geometry is related to but does not cover the geometry of Liouville vector fields.
   Below we discuss the case of Liouville vector fields coming from regular sections of the prequantization bundle only.

   The core $N(\lambda)$ for a given Liouville vector field $\lambda$ consists of finite trajectories of the flow $\Phi^t_{\lambda}$.
   Thus the core $N(\lambda_{\alpha}) \subset M \backslash D_{\alpha}$ defined by a regular section $\alpha$,  can be characterized by the following properties: at smooth point $p \in N(\lambda_{\alpha})$ vector field $\lambda_{\alpha}$ is tangent to $N(\lambda_{\alpha})$; the core $N(\lambda_{\alpha})$ is stable with respect to the flow $\Phi^t_{\lambda_{\alpha}}$; the core $N(\lambda_{\alpha})$ does not touch $D_{\alpha}$.

   Therefore for Lioville vector field $\lambda_{\alpha}$, defined by a regular section $\alpha \in \Gamma (M, L)$ one has the following observation:
   if $S$ is a SBS lagrangian submanifold with respect to $[\alpha]$ then $S$ is contained by $N(\lambda_{\alpha})$. Indeed, as we have mentioned above $\lambda_{\alpha}$ is tangent
   to $S$ if it is SBS with respect to $[\alpha]$, therefore $S$ must be stable with respect to the flow; on the other hand by the very definition $S$ does not touch $D_{\alpha}$.
   In the opposite direction we can just say that if the core $N(\lambda_{\alpha})$ contains a $n$ - dimensional smooth component $N_i$ which is lagrangian then $N_i$
   must be SBS.

   However it is known that $N(\lambda)$ can have bigger dimension than $n$. The desired bound takes place if a given Liouville vector field $\lambda$ is gradient like for a smooth function $\phi$ namely it exists a compatible rimannian metric $g$ such that 
   $$
   \rm{d} \phi(\lambda) \geq C \Vert \lambda \Vert^2_g
   \eqno (2)
   $$
    for certain positive constant $C> 0$. One says that such a pair
   $(\phi, \lambda)$ defines Weinstein structure on the open symplectic manifold $M \backslash D$. In this case the core $N(\lambda)$ is called Weinstein skeleton $W(M \backslash D)$:
   it is at most $n$ - dimensional and isotropic. Therefore one has

   {\bf Proposition 9.} {\it If $\lambda_{\alpha}$ admits a smooth function $\phi$ such that $(\phi, \lambda_{\alpha})$ induces a Weinstein structure on $M \backslash D_{\alpha}$ then
   a smooth $n$ - dimensional submanifold $S$ is SBS with respect to $[\alpha]$ if and only if $S$ is a component of Weinstein skeleton $W(M \backslash D_{\alpha})$.}

   The basic examples of Weinstein structures come from the complex geometry: if it exists a compatible almost complex structure $I$ then every pseudo holomorphic
   section $\alpha \in \Gamma (M, L)$ induces a Weinstein structure on the complement $M \backslash D_{\alpha}$: the function $\phi$ equals to $\psi_{\alpha} =
   - \rm{ln} \vert \alpha \vert_h$, and in this case the Liouville vector field $\lambda_{\alpha}$ is the gradient vector field for this $\phi$
   with respect to Riemannian metric $g$, reconstructed from $\omega$ and $I$. However as we have mentioned above for integrable $I$ the corresponding
   Weinstein skeleton is not smooth: it is just a CW complex, which does not admit smooth closed components. At the same time
   every closed $n$ - dimensional cycle in $W(M \backslash D_{\alpha})$ presents a non trivial homology class in $H_n(M \backslash D_{\alpha}, \mathbb{Z})$
   due to the dimensional reason. 
   
   Our main interest in the story above was the following: is it possible to find a version of {\bf $\mathbb{P}$ - covering theorem}, so is it possible to
   reconstruct a family of pairs $([\alpha]_t, S_t) \in {\cal M}_{SBS}(a_I)$ starting with a path $\{ [\alpha_t] \} \subset \mathbb{P}(D_{\alpha})$
   (as it was shown above it is essentially the same as for a path in whole $\Gamma (M, L)$)? In certain cases the answer is definitely yes:
   for example if the path $\{ [\alpha_t \}$ is given by a Hamiltonian deformation of a finite region in $M \backslash D_{\alpha}$ since for this case
   we can deform the pair $([\alpha_1], S_1)$ exploiting Proposition 6. 
   
   Another particular case is given by the variations of Weinstein structures: suppose that  Liouville vector field $\lambda_{\alpha_1}$ is given  by
   a section $\alpha_1$  which is pseudo holomorphic with respect to an almost complex structure $I$ on the complement $M \backslash D_{\alpha_1}$, therefore
    if we take $\phi = \psi_{\alpha_1}$ then the pair $(\phi, \lambda_{\alpha_1})$
   defines a Weinstein structure. Then we can vary the second element in the pair such that the new pair again defines a Weinstein structure.
   For example consider the following variation $\lambda_t = \lambda_{\alpha_1} + (t-1) X_{\phi}, t \in [1; + \infty)$ where $X_{\phi}$ is the Hamiltonian vector field
   of function $\phi$ (moreover we can take any composite function derived from $\phi$): it is not hard to see that
   for any $t \geq 1$ it exists a positive constant $C_t$ such that the pair $(\phi, \lambda_t)$ satisfies condition (2) with respect to the same metric.
   Indeed,  $\rm{d} \phi(\lambda_t) \equiv \rm{d} (\lambda_1)$ and $\vert \lambda_t \vert_g^2 = (1 + t^2) \vert \lambda_1 \vert^2_g$
   since $X_{\phi}$ and $\lambda_1$ are related as Hamiltonian vector field and gradient vector field for the same function. Therefore
   $C_t$ can be taken equal to  $C(1+t^2)$.     

   Geometrically this means that we deform components of the Weinstein skeleton along the level sets of the potential $\phi$, fixing the critical points
   which are stable under the process. 
   
   On the other hand the first necessary property which follows from (2) says that $\lambda$ must have more critical points than in general case
   since if $x \in M \backslash D$ is a critical point of $\phi$ which is not a zero point for $\lambda$ we immediately get a contradiction.
   Therefore $\{ \rm{Crit} \phi \} \subseteqq \{ \rm{sing} \lambda \}$. Note that if we fix any set of points $\{p_1, ..., p_l \} \subset M \backslash D$
   we can define an ideal in $C^{\infty}(M, \mathbb{R})$ consists of smooth functions  $F$ such that $\rm{Crit} F \supseteqq \{p_1, ..., p_l \}$.
   Indeed, if $\rm{d} F_1 (p_i) = \rm{d} F_2 (p_i) = 0 $ it follows $\rm{d}(F_1 + c F_2)(p_i) = 0$ and $\rm{d} (F_1 \cdot F_2) (p) = 0$.
   Then for a given Weinstein structure $(\phi, \lambda)$ one has the following variation space: take ideal ${\cal I}(\lambda) \subset C^{\infty}(M, \mathbb{R})$
   taking the singular points $p_1, ..., p_l$ of the vector field then it exists a small neighborhood of zero such that for any small function $\delta
   \in {\cal I}(\lambda)$ the pair $(\phi + \delta, \lambda)$ again defines a Weinstein structure. 
   
   But again the critical points for the Weinstein skeleton remain to be the same under such deformation, therefore all the types
   can help just in the testing of the existence problem for Bohr - Sommerfeld resolutions.
   
   To find a way how  a global version of the covering theorem can be formulated we study the following situation, rather close
   to the cases, considered above. Let $[\alpha]$ be the class of holomorphic section, such that the corresponding Weinstein skeleton
   $W(X \backslash D_{\alpha})$ admits $n$ - dimensional closed cycles. Let $S_1$ be a smooth lagrangian submanifold in $M \backslash D_{\alpha}$
   such that it is SBS for some tramsformation $P(F) [\alpha]$. Can we join $[\alpha]$ and $P(F) [\alpha]$ in $\mathbb{P}(D_{\alpha})$
   by a path $\{[\alpha] + \rho_t \}$ in $\rho$ - representation such that for each $t$ one has  SBS submanifold $S_t$ which
   satisfies $(\rm{Im} \rho(\alpha) + \rho_t)|_{S_t} \equiv 0$? 
   
   First of all note that $S_1$ is necessary {\it exact} with respect to the Weinstein structure, defined by the holomorphic section. Indeed,
   the restriction $\rm{Im} \rho(\alpha)|_{S_1} = - \rm{d} F|_{S_1}$. As the first step we can vary the function $F$ such that its value on $S_1$
   is unchanged but the restriction of the corresponding vector field $\lambda_1 = \omega^{-1}(\rm{Im} \rho(\alpha) + \rm{d} F)$ would
   be suitable for the further analysis. Namely the function $F$ can be chosen such that the corresponding vector field $\lambda_1$ 
   {\bf} identically vanishes on $S_1$. Indeed, take a function $F_0$ such that the restriction $(\rm{Im} \rho(\alpha) + \rm{d} F_0)_{S_1}$ is trivial, then
   the corresponding vector field $\omega^{-1}(\rm{Im} \rho(\alpha) + \rm{d} F_0)$ must be parallel to $S_1$ at its points. For a Darboux - Weinstein
   neighborhood ${\cal O}_{DW}(S_1)$ the pair (function $f$ on $S_1$, vector field $v$ on $S_1$) generates a smooth function  $F(x, p)
   = f(x) + p(v)$; it is not hard to see that for our situation the function $F(x, p)$, constructed for $f = F_0|_{S_1}$ and
   $v = \omega^{-1}(\rm{Im} \rho(\alpha) + \rm{d} F_0)|_{S_1}$,  posses the property $X_{F} \equiv - \lambda$ at the points of $S_1$
   where $\lambda = \omega^{-1}(\rm{Im} \rho(\alpha)$ --- the second element, defined the Weinstein structure. 
   
   Therefore we can find the transformation term $F$ such that $\rm{Im} \rho(\alpha) + \rm{d} F$ coincides to the canonical 1 - form $\alpha_{can}$ being
   restricted to the Darboux - Weinstein neighborhood ${\cal O}_{DW}(S_1)$ (essentially it repeats the arguments we have used in the proof of Proposition 5 above).
   
   Given Weinstein structure $(\psi_{\alpha}, \lambda)$ induces certain numerical attachment for the space of compact exact lagrangian submanifolds:
   since the restriction of function $F|_{S_1}$ is unique then one can take the maximal and minimal values of this restriction and define
   $$
   N(S_1) = \rm{max} F|_{S_1} - \rm{min} F_{S_1}.
   $$
   Note that this number does not depend on any choice rather the Weinstein structure itself and a given exact lagrangian submanifold. 
   
   Note that this number is great or equal to 0, and the last can happen if and only if $S_1$ is a component of the Weinstein skeleton $W(M \backslash D_{\alpha}$.
   One can understand this number as a {\bf distance} from $S_1$ to $W(M \backslash D_{\alpha}$.
   
   Now our aim is the find a procedure which starts with a given $S_1$ and results with a smooth exact lagrangian submanifold $S_{T}, T \in [0; 1)$, such that
   $$
   N(S_1) \cdot T = N(S_T).
   $$
   
   Consider the Darboux - Weinstein neighborhood ${\cal O}_{DW}(S_1)$. Recall, see [5], that for every small Bohr - Sommerfeld deformation of $S_1$
   corresponds to a smooth function $f \in C^{\infty}(S_1, \mathbb{R})$ by the following rules: if the deformed $\delta S$ belongs to the DW - neighborhood,
   then it is defined by the graph of an exact 1 - form on a small neighborhood of zero section in $T^* S_1$ therefore it is defined by a smooth function
   $f \in C^{\infty}(S_1, \mathbb{R})$ up to constant. Note however that the Darboux - Weinstein neighborhood is not unique, but we have fixed appropriate one above.
   All small Bohr - Sommerfeld deformations give us exact lagrangian submanifolds, therefore we get a correspondence 
   $$
   \Lambda_{\alpha}: {\cal O}_{\varepsilon} (C^{\infty}(S_1, \mathbb{R}) \to C^{\infty}(S_1, \mathbb{R})
   \eqno (3)
   $$
   from $\varepsilon$ - neighborhood of constant functions with respect to the norm $\vert f \vert^{max}_{min} = \rm{max} f - \rm{min} f$ on $S_1$. Namely for a function
   with "small differential" $f$ we take the corresponding Bohr - Sommerfeld deformation $S_f \subset {\cal O}_{DW}(S_1)$, restrict the form $\rm{Im} \rho(\alpha)$
   to $S_f$ and by the definition get an exact 1 -form on $S_f$. This 1 - form is the differential of a smooth function $\phi$ on $S_f$, which can be
   chosen say  to have the same maximal value as $f$. Then using the Darboux - Weinstein neighborhood structure we can project $S_f$ to $S_1$ and get
   the corresponding function $\Lambda_{\alpha} (f)$ on $S_1$.      
   
   It is not hard to see that $F|_{S_1} = - \Lambda_{\alpha}(0)$, 
   $$
   N(S_1) = \vert F|_{S_1} \vert^{max}_{min},
   $$
   and for a Bohr - Sommerfeld deformation $S_f$ one has $N(S_f) = \vert \Lambda_{\alpha} (f) \vert^{max}_{min}$ by the very definition. Therefore to find desired $S_T$, described
   above, we need to find a function with small differential $\delta \in C^{\infty}(S_1, \mathbb{R})$ such that $\vert \Lambda_{\alpha} (\delta) \vert^{max}_{min} <
   \vert \Lambda_{\alpha} (0) \vert^{max}_{min}$. Therefore we combine the map (3) with the max - min norm of the images and define
   $$
   N_{\alpha}: {\cal O}_{\varepsilon} C^{\infty} (S_1, \mathbb{R}) \to \mathbb{R}_{\geq 0}.
   $$
   
   For this map we have the following:
   
   {\bf Proposition 10.} {\it The map $N_{\alpha}$ is smooth. The only  possible critical value is zero.}
   
   This Proposition implies the existence of deformation $S_T$ of a given exact lagrangian submanifold $S_1$, which decreases the distance to the Weinstein skeleton:
   since the only critical value is zero we can not decrease the distance if and only if the distance is already zero so $S_1$ itself belongs to $W(M \backslash D_{\alpha}$.
   Note that $S_T$ is  Hamilton isotopic to $S_1$ therefore if $S_T$ already {\bf regular} to implies that $S_1$ is regular too.

   Geometrical essence of the picture can be described as follows. Consider the moduli space ${\cal B}_S$ of Bohr - Sommerfeld lagrangian submanifolds constructed in
   [5]. Suppose we have a Weinstein structure $(\phi, \lambda)$ on the complement $M \backslash D$ where $D$ is a symplectic submanifold.
    Then it is defined locus $\Delta \subset {\cal B}_S$
   consists of such $S$ that the intersection $S \cap D$ is non empty. Then the complement ${\cal B}_S \backslash \Delta$ is a set of connected components
   ${\cal B}_1, ..., {\cal B}_k, ...$. Suppose that a component ${\cal B}_i$ contains an exact lagrangian submanifold $S_1$, then it implies that every $S$
   from the same component is exact too. Therefore one has a vector field $\Theta_{\cal B}(\lambda) \in \rm{Vect} {\cal B}_i$ defined by exact 1 - forms
   $\omega(\lambda)|_{S}$ since each exact 1 - form on $S$ corresponds to tangent vector, see [5]. Max - min norm gives a normalization
   of the vector field $\Theta_{\cal B}(\lambda)$. Note that when $S$ tends to locus $\Delta$ the norm of $\Theta_{\cal B}(\lambda)([S])$ goes to infinity:
   the main reason is that $S$ can not collapse being Bohr - Sommerfeld and can not  in the limit be contained by $D$ therefore the maximal value
   shall go to infinity but the minimal must remain finite. Proposition 10 says that if for $[S_1] \in {\cal B}_i$ the tangent vector
   $\Theta_{\cal B}([S_1])$ is non trivial then in a small neighborhood of $[S_1]$ in ${\cal B}_i$ it exists some smooth Bohr - Sommerfeld submanifold
   $S_T$ such that the norm of the tangent vector $\Theta_{\cal B}(\lambda)([S_T])$ strictly less than for $[S_1]$. Due to the observation above
   this deformation gives submanifolds which are  even further away from $D$ than $S_1$ so this deformation does not affect the basic property we have exploited in the construction.
   
    To prove  Proposition 10 we can simplify the picture as follows: consider cotangent bundle $T^*S$ of a smooth $n$ - dimensional compact connected manifold $S$ together with the canonical  action 1 - form $\alpha_{can}$ such that the canonical symplectic form is given by $\omega_{can} = \rm{d} \alpha_{can}$. For a given global smooth function $F(x, p)
    \in C^{\infty}(T^* S, \mathbb{R})$     one has   a map $\Lambda_F: C^{\infty}(S, \mathbb{R}) \to C^{\infty}(S, \mathbb{R})$ defined  by the procedure:
  restrict $F$ to the graph $\Gamma (\rm{d} f) \subset T^*S$ then project the restriction $F_{\Gamma(\rm{d} f)}$ to $S$ using the canonical projection $\pi: T^*S \to S$ 
  and at the end add to the result the function $f$ itself. Due to the Darboux - Weinstein theorem all facts about $\Lambda_{\alpha}$ above follow from the  same
   statements for $\Lambda_F$ (note that the restriction $\alpha_{can}|_{\Gamma(\rm{d} f)}$ equals to $\pi^*(\rm{d} f)$ by the very definition), therefore we 
   can study the last map, establish certain properties and then extend the results to the situation of Proposition 10. 
   
   The first property of map $\Lambda_F$ can be easily established: the map is injective. Indeed, suppose that for two functions $f_1$ and $f_2$ one gets
   $\Lambda_F(f_1) \equiv \Lambda_F(f_2) \in C^{\infty}(S, \mathbb{R})$. First, suppose that $f_2 = f_1 + const$, then the graph $\Gamma(\rm{d} f_i)$
   is the same therefore the restriction $F|_{\Gamma(\rm{d} f_i)}$ is the same, consequently it can happen if and only if $f_2 - f_1 \equiv 0$.
   
   Now in the other case  $f_1$ and $f_2$ define two distinct graphs $\Gamma(\rm{d} f_i)$ which are lagrangian submanifolds in $T^*S$ and since $S$ is compact
   they must intersect each other at least at two points $p_+, p_-$. It follows from the fact that non constant function $f_2 - f_1$ must have
   at least two critical points on the compact manifold $S$ --- maximal and minimal, and these points underly the intersection points $p_+, p_-
   \in \Gamma(\rm{d} f_2) \cap \Gamma(\rm{d} f_1)$. Since $S$ is connected we can choose two pathes $\gamma_i \subset \Gamma(\rm{d} f_i)$
   with ends at $p_{\pm}$. Since one supposes $\Lambda_F(f_1) = \Lambda_F(f_2)$ it follows that the integrals 
   $$
   \int_{\gamma_i} (\alpha_{can} + \rm{d} F)|_{\Gamma(\rm{d} f_i)}
   $$
   must be the same; therefore $f_2 (p_{\pm}) = f_1 (p_{\pm})$. But the points were chosen as maximal and minimal for the difference $f_2 - f_1$ which implies
   $f_2 - f_1 \equiv 0$. 
   
   In the finite dimensional case we shall have the inverse map $\Lambda_F^{-1}$ which would solve;  in the present case we are interested only in the decreasing
  of the $\vert^{max}_{min}$ - norm  for the restriction of $F$ to the zero section. It is clear that
   if $F|_{S} = 0$ or other constant then $\Lambda_F(0)$ equals to the same constant therefore the image $\Lambda_F (f)$ corresponds to any deformation $\Gamma(\rm{d} f)$
   must be non constant due to the injectivity property therefore $\Lambda_F(f)$ must have  non zero $\vert^{max}_{min}$ - norm. At the same time
   if $F$ being restricted to the zero section is non constant then it exists a deformation $\delta f$ such that $\vert \Lambda_F(0) \vert^{max}_{min} > \vert \Lambda_F(\delta f) \vert^{max}_{min}$.
   
   To establish this fact consider infinitesimal variation $\delta f \in C^{\infty}(S, \mathbb{R})$ and the corresponding deformation $\Gamma(\rm{d}\delta f)$ of the zero section in
   $T^* S$. Then the restriction $F|_{\Gamma(\rm{d} \delta f)}$ has the linear part $F(x, 0) + \rm{d} \delta f (Y_F)$ where $Y_F$ is a vector field on the zero section
   given by $\sum_{i=1}^n \frac{\partial F}{\partial p_i} \frac{\partial}{\partial x_i}$ in a local Darboux  coordinate system. Therefore the differential
   of the map $\Lambda_F$ at zero is presented by the linear map $\delta f \mapsto \rm{d} \delta f (Y_F) + \delta f$. Note that the differential admits trivial kernel:
   indeed, if for certain non constant $\delta f$ the result is constant then at the maximal and minimal points of $\delta f$ the result must be equal to both the values
   therefore $\delta f$ has to be constant. 
   
   If $F(x, 0)$ is non constant then we can exploit the following strategy to find a deformation $\delta f$
  such that the norm $\vert \Lambda_F(\delta f) \vert^{max}_{min}$ would be less than the norm of the initial $F(x, 0)$. Suppose that maximal and minimal points
  $p_{\pm}$ of  $F(x, 0)$ are isolated, then choose small neighborhoods of $p_{\pm}$ and take bump functions which are close to $- t F(x, 0)$ in the neighborhoods
  and which are trivial outside of the neighborhoods where $t$ is a small parameter. Then the sum of such bump functions has zero differential at $p_{\pm}$ and
  therefore the part $\rm{d} \delta f(Y_F)$ is killed at the points; and the same time the maximal and minimal values of $\Lambda_F$ for such $\delta f$ must be
  less than for our given $F(x, 0)$. 
  
  These arguments leads to the proof of  Proposition 10; however we can not present an appropriate answer to the following
  
  {\bf Problem.} {\it Find an effective bound on possible  decreasing of the norm $\vert \Lambda_{\alpha} (\delta f) \vert^{max}_{min}$ for a given
  function $F$.}
  
  The dreamed  form for the bound: it is a constant $C$ which depends on the topology of a given exact lagrangian submanifold $S$ and the properties of  $\alpha$
  such that  it exists a deformation $S_T$ of $S$ satisfying $\vert \rm{Im} \rho(\alpha)|_{S_T} \vert^{max}_{min} \leq C \vert \rm{Im} \rho(\alpha)|_S \vert^{max}_{min}$.
  
  If the problem admits such a solution then it were possible to prove  a version of $\mathbb{P}$ - covering theorem, mentioned above. Moreover,
  it were possible to attack Eliashberg conjectures, see [8], for the case when Weinstein structures are given by good smooth sections of
  the prequantization bundle (in particular for pseudo holomorphic sections). Indeed, as we have discussed above each exact lagrangian submanifold
  $S$ represents a point in the moduli space ${\cal B}_S$, and the corresponding form $\rm{Im} \rho(\alpha)$ is presented in the form
  $\alpha_{can} + i \rm{d} F$ being restricted to a Darboux - Weinstein neighborhood of $S$. If in the neighborhood one can find a deformation $S_T$
  such that the norm of the restriction $\rm{Im} \rho(\alpha)|_{S_T}$ subjects our "dreamed bound" then it is possible to find a chain
  of exact lagrangian submanifolds $S, S_{T_1}, S_{T_2}, ...$ such that the norms of the restrictions $\rm{Im} \rho(\alpha)|_{S_{T_i}}$ tends to zero,
  therefore $S_{T_i}$ is closer and closer to the Weinstein skeleton $W(M \backslash D_{\alpha})$. In the limit $\{ S_{T_i} \}$ touches the skeleton
  therefore it must be the corresponding limiting cycle in the CW - complex, and this cycle must be non trivial by the dimensional reasons.
  It would imply that the starting submanifold $S$ must present non trivial homology class in $H_n(M \backslash D_{\alpha}, \mathbb{Z})$.
  Furthermore, using the deformations of the Weinstein structures discussed above one can try to catch certain slight  deformation such that
  $S_{T_i}$ being sufficiently closed to $W(M \backslash D_{\alpha})$ is regular with respect to this slightly deformed Weinstein structure. But
  all $S_{T_i}$ are Hamiltonian isotopic on the complement $M \backslash D_{\alpha}$ by the construction therefore
  this slightly deformed Weinstein structure can be transport back to $S$, which would lead to the fact that $S$ is regular itself.
  
  $$$$
  
  $$$$
  
  {\bf References:}
  
  [1] N.A. Tyurin, {\it "Special Bohr–Sommerfeld Lagrangian submanifolds"},  Izv. Math., 80:6 (2016), 1257–1274; 
  
  [2] N. A. Tyurin, {\it "Special Bohr–Sommerfeld Lagrangian submanifolds of algebraic varieties"},  Izv. Math., 82:3 (2018), 612–631;
  
  [3] N.A. Tyurin, {\it "The moduli space of D-exact Lagrangian submanifolds"},  Siberian Math. J., 60:4 (2019), 709–719; 
  
  [4] N.A. Tyurin, {\it "On the Kahlerization of the Moduli Space of Bohr–Sommerfeld Lagrangian Submanifolds"},  Math. Notes, 107:6 (2020), 1008–1009;
  
  [5] A.L. Gorodentsev, A.N. Tyurin, {\it "Abelian Lagrangian algebraic geometry"},  Izv. Math., 65:3 (2001), 437–467;  
  
  [6] K. Cieliebak, Ya. Eliashberg, {\it "From Stein to Weinstein and back - symplectic geometry of affine complex geometry"},
   Colloq. Publ., 59,  Amer. Math. Soc., Providence, RI, 2012;
  
  [7]  N. A. Tyurin, {\it "Dynamical correspondence in algebraic Lagrangian geometry"}, Izv. Math., 66:3 (2002), 611–629;
  
  [8]  Ya. Eliashberg, {\it "Weinstein manifolds revisited"},  Modern geometry: a celebration of the work of Simon Donaldson, 59–82,
   Proc. Sympos. Pure Math., 99, Amer. Math. Soc., Providence, RI, 2018.

\end{document}